\documentclass[conference]{IEEEtran}
\usepackage{multirow}
\usepackage{multicol}
\usepackage{lipsum}
\usepackage{graphicx}
\usepackage{graphics}
\usepackage{amsmath}
\usepackage{amsbsy}
\usepackage{blindtext}
\usepackage{tcolorbox}
\usepackage{amssymb}
\usepackage{scalerel}
\usepackage{mathabx}
\usepackage{verbatim}
\usepackage{booktabs}
\usepackage{rotating,tabularx}
\usepackage{mathrsfs} 
\usepackage{hyperref}

\usepackage{graphicx}
\usepackage{epstopdf}
\usepackage{algorithm}

\usepackage{algpseudocode}

\usepackage{adjustbox}
\usepackage{soul}
\usepackage{xcolor}
\usepackage{cite}
\usepackage{tikz}
\usepackage{color, colortbl}
\usepackage[first=0,last=9]{lcg}
\definecolor{LightGray}{rgb}{0.7,0.7,0.7}

\usepackage[wby]{callouts}
\usetikzlibrary{shapes.multipart}

\usepackage{array}
\usepackage{makecell}

\allowdisplaybreaks

\usepackage{amsthm}
\theoremstyle{definition}

\theoremstyle{remark}

\usepackage[utf8]{inputenc}
\usepackage[english]{babel}

\usepackage[caption=false,font=footnotesize]{subfig}
\usepackage[T1]{fontenc}
\usepackage{scalerel,stackengine}
\newcommand\reallywidecheck[1]{%
\savestack{\tmpbox}{\stretchto{%
  \scaleto{%
    \scalerel*[\widthof{\ensuremath{#1}}]{\kern-.6pt\bigwedge\kern-.6pt}%
    {\rule[-\textheight/2]{1ex}{\textheight}}
  }{\textheight}%
}{0.5ex}}%
\stackon[1pt]{#1}{\scalebox{-1}{\tmpbox}}%
}
\stackMath
\IEEEoverridecommandlockouts
\IEEEoverridecommandlockouts

\definecolor{shadecolor}{RGB}{190,190,190}

\renewcommand{\baselinestretch}{1.05}

\newif\ifarxiv
\arxivtrue
\arxivfalse

\begin{document}

\title{\LARGE\bf Demand Charge Management: Prototype Design and Testing}

\author{Christian Escarrega, Wyatt Lopez Coggins, Sammy Hamed, \\Mohammadreza Iranpour, Mohammad Rasoul Narimani


\thanks{This work has been supported from NSF contract \#2308498.

979-8-3315-4112-5/25/\$31.00 ©2025 IEEE}%
}

\maketitle

\begin{abstract}
This paper presents the design, implementation, and validation of a smart, low-cost Energy Management System (EMS) and Demand Charge Management (DCM) prototype, developed as part of an undergraduate senior design project. The system serves as both a practical solution for reducing electricity costs and a pedagogical tool for teaching real-time energy control concepts in power and embedded systems courses. Unlike conventional EMS/DCM solutions that rely on high-cost commercial hardware or purely theoretical models, the proposed system integrates grid power, lithium-iron phosphate (LiFePO$_4$) battery storage, and real-time control into a unified, scalable platform constructed at a fraction of the cost, approximately \$1,800 compared to over \$16,000 for leading commercial options. The controller dynamically optimizes energy usage by switching between grid and battery sources based on real-time measurements of electricity prices, load power, and battery state of charge (SoC). This enables peak shaving, energy arbitrage, and backup power functionality, thereby enhancing cost efficiency and grid resilience for both residential and small commercial users. The architecture features a modular three-layer design comprising a sensing layer for electrical data acquisition, a control layer executing Python-based logic on a Raspberry Pi, and an actuator layer for seamless energy switching. Data is communicated via MQTT and visualized through the Blynk IoT platform, providing an intuitive and remotely accessible user interface. The prototype's effectiveness was validated through real-world testing, confirming its capability to reduce demand charges and ensure reliable energy delivery under varying operational conditions. Its affordability, open-source control logic, and educational versatility make it an ideal candidate for both deployment and instructional use.

\end{abstract}


\section{Introduction}
\label{Introduction}

The increasing complexity of modern power distribution systems and the steady rise in electricity costs have made effective electricity demand management more critical than ever. Residential and commercial consumers often face significant financial pressure due to high electricity bills, especially during peak demand months like summer, when air conditioning usage surges \cite{Zhao2017,Li2019}. In response, utilities employ various pricing strategies, such as time-of-use (TOU) tariffs and demand charges, to regulate consumption and stabilize the grid \cite{Al-Wakeel2017}.
These strategies typically include several components. TOU pricing applies different rates based on time-of-day, with higher charges during peak hours \cite{Xu2020}. Additional charges may apply when users exceed predefined energy thresholds \cite{Bitar2017}. Seasonal price fluctuations further complicate billing, as energy costs rise during periods of extreme heating or cooling needs \cite{Ma2022, asghari2019optimal,narimani2018optimal,narimani2017energy, narimani2015effect,niknam2012efficient, narimani2015dynamic, narimani2018optimal, azizivahed2017new, narimani2016reliability}. These pricing mechanisms, while intended to encourage responsible energy use, often result in unpredictable and burdensome energy bills for consumers.

To address these challenges, Demand Charge Management (DCM) systems have emerged as promising tools. These systems help users optimize energy consumption by intelligently balancing power draw between grid and off-grid sources, such as batteries or local generators \cite{Wang2020,Chauhan2018}. A well-designed DCM controller can not only reduce peak-hour electricity expenses but also improve energy efficiency and reliability by shifting load to off-peak periods \cite{Paudyal2011}. Moreover, such systems can support grid stability and contribute to broader sustainability goals by encouraging the integration of renewable energy sources \cite{Huang2019,Mohamed2020}. For example, Narimani \cite{narimani2019demand} proposed a residential energy management strategy that minimizes cost by rescheduling appliance usage in response to hourly price signals, while maintaining user comfort. In another study, Asghari et al. \cite{asghari2019method} introduced a method for controlling distributed energy storage systems based on multi-objective optimization that balances demand charge reduction with photovoltaic (PV) utilization.

Recent developments have further refined these concepts by incorporating user behavior analytics and 24-hour operational cycles. For instance, intelligent controllers have been designed to autonomously discharge batteries when consumption nears peak thresholds and recharge during low-demand periods, enhancing both cost savings and operational efficiency \cite{Tang2021}. Various advanced control techniques, including model predictive control (MPC) and fuzzy logic, have also been proposed to optimize energy consumption in residential and industrial settings \cite{Zhou2022,Beaudin2015}. Complementing these are Internet-of-Things (IoT)-based solutions that employ real-time monitoring and machine learning to adaptively manage energy use \cite{Lin2018}.

Despite these advances, many existing systems remain either simulation-based or prohibitively expensive for average consumers. They often lack practical, hardware-based implementations that simultaneously address economic constraints and real-time demand control \cite{Chen2020}. This gap highlights the need for a cost-effective, hardware-driven DCM solution that can be adopted in both residential settings and educational environments.
In this paper, we propose a smart, low-cost Energy Management System (EMS) and DCM prototype developed by undergraduate students as part of a senior design project. Unlike high-cost commercial systems such as the Tesla Powerwall 3, Enphase IQ Battery 5P, and Generac PWRcell, which often exceed \$16{,}000 \cite{TeslaPowerwall3}, our system was built at a component cost of approximately \$1{,}800. The cost reduction makes the platform accessible and ideal for teaching power and embedded systems.

\begin{figure*}
    \centering
\includegraphics[scale=0.34,trim=2cm 8.5cm 0.1cm 9cm,clip]{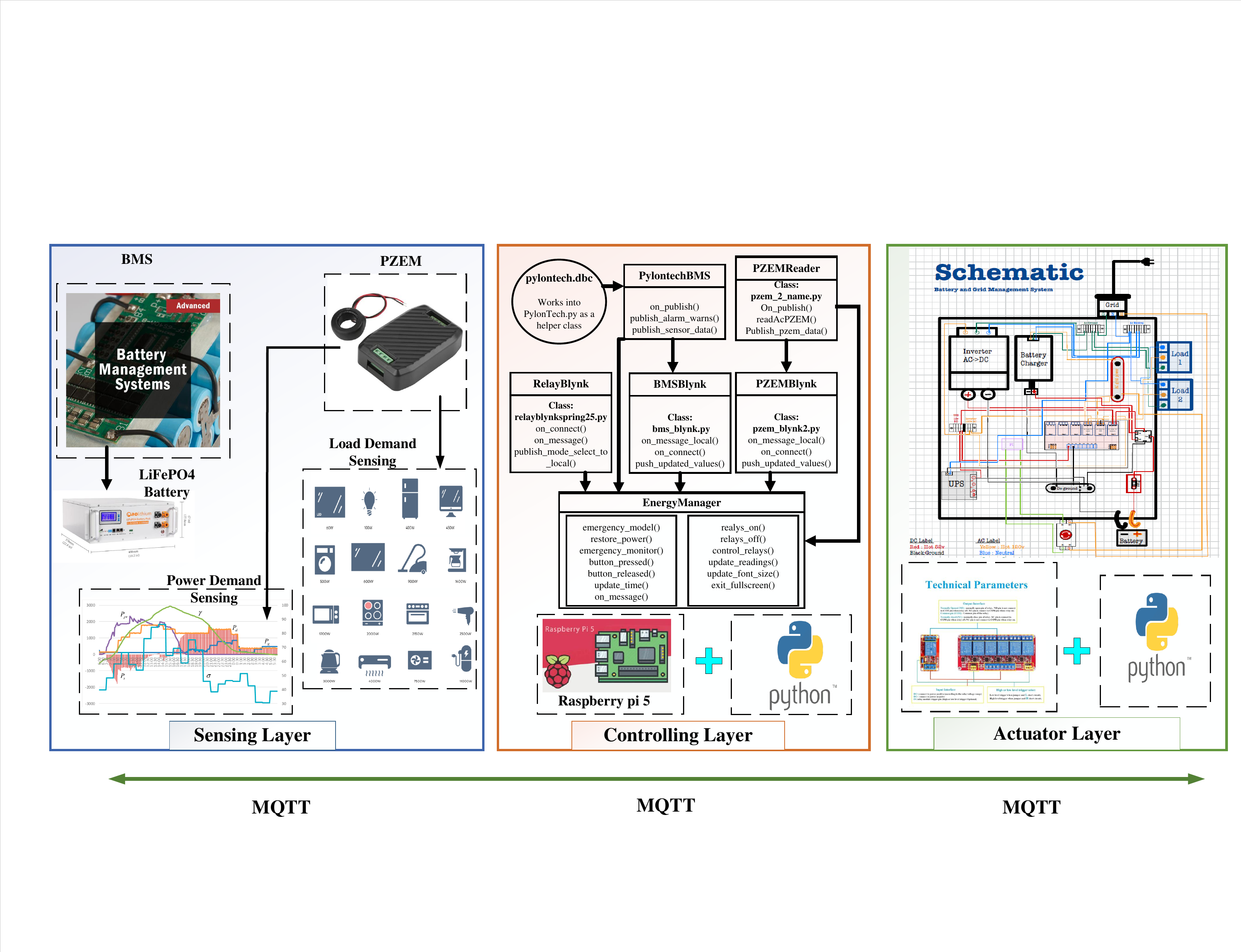}
    \caption{Schematic of the proposed framework for Demand Charge Managemen.}
    \label{fig:mainfigure}
\end{figure*}

The proposed system intelligently integrates grid power, local generation, and battery storage to minimize electricity costs through two main control parameters: (1) the battery's State of Charge (SoC), and (2) predefined power thresholds for peak demand control. Energy usage is dynamically optimized by switching between power sources based on real-time data related to load demand, electricity pricing, and SoC.
To implement this solution, we developed a three-layer controller architecture consisting of a sensing layer for electrical data acquisition, a control layer executing Python-based decision logic on a Raspberry Pi, and an actuator layer that performs energy source switching. The system is further integrated with MQTT communication and the Blynk IoT platform for real-time data visualization and user interaction.
This paper details the design, implementation, and experimental validation of the prototype. The results demonstrate the system's effectiveness in demand charge mitigation, energy cost reduction, and educational applicability.

\section{Methodological Framework}
\label{sec:Methodological Framework}

This section outlines the design and implementation of the proposed DCM system, emphasizing its potential to enhance energy efficiency and cost-effectiveness for both residential and commercial users. The system is engineered to intelligently regulate power flow among the utility grid, battery storage, and end-use loads, thereby minimizing electricity expenses while maintaining reliable service.
To realize this functionality, a three-layer controller architecture was developed. This architecture integrates real-time monitoring, intelligent decision-making, and automated energy dispatch into a unified, hardware-based framework. The three functional layers, depicted in Fig.~\ref{fig:mainfigure}, are described as follows:

\textbf{Sensing Layer:}  
The sensing layer is responsible for continuous monitoring of electrical parameters across the system. It captures key metrics such as voltage, current, active power, and power factor from both the supply side (grid and battery) and the demand side (load). The battery's SoC is also tracked in real-time using a Battery Management System (BMS). These data points serve as essential inputs for the control logic, enabling precise and up-to-date decision-making.

\textbf{Control Layer:}  
At the core of the system is the control layer, implemented using a Raspberry Pi 5. This layer executes custom Python-based algorithms to analyze sensor data and determine the optimal power source, either grid or battery, based on predefined thresholds and time-of-use considerations. The control logic dynamically manages power dispatch by evaluating SoC, load demand, electricity pricing, and user-configured modes. This enables the system to balance energy cost, battery health, and grid reliability while allowing for both autonomous and user-override operations.

\textbf{Actuator Layer:}  
The actuator layer physically executes the control decisions by regulating energy flow through a set of relays, converters, and switching circuits. It enables seamless transitions between power sources, mitigates supply disruptions, and enforces load-shifting strategies. For instance, when load demand exceeds a predefined threshold during peak pricing periods, the system prioritizes battery discharge. Conversely, it switches back to grid power during off-peak periods or when battery SoC is critically low. The actuator layer also integrates emergency stop functionality to ensure operational safety.

By integrating these three layers into a cohesive hardware prototype, the proposed system offers a cost-effective, intelligent, and scalable solution for real-time energy management. It provides a practical platform for reducing electricity costs, enhancing energy resilience, and supporting instructional use in power systems and embedded control education.

\section{High-Level Fidelity Model}
\label{sec:high level model}

\begin{figure}
    \centering
\includegraphics[scale=0.34,trim=27.5cm 3cm 10cm 1cm,clip]{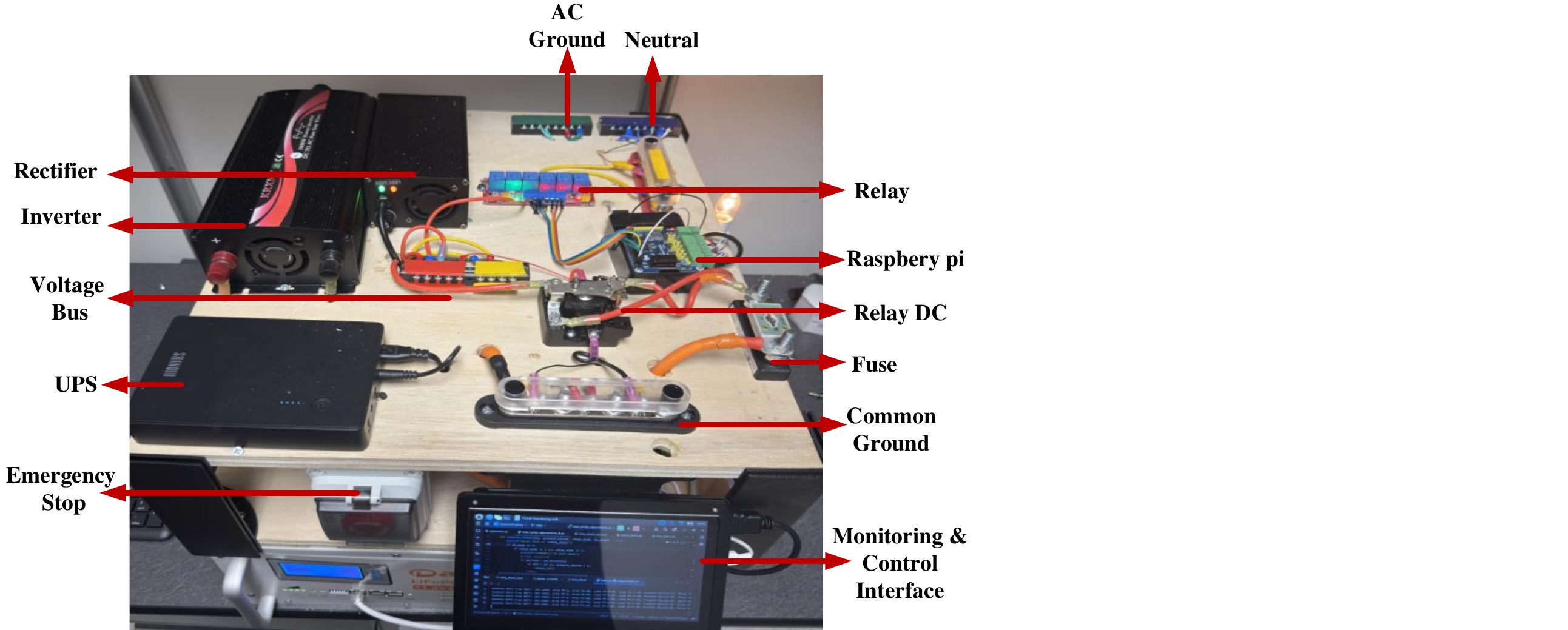}
    \caption{Prototype Demand Charge Managemen system.}
    \label{fig:prototype}
\end{figure}

The proposed Demand Charge Management (DCM) system integrates a suite of hardware and software components designed to optimize electricity usage and minimize peak demand costs. As previously outlined, the system intelligently manages energy flow between the utility grid and a local battery storage unit through coordinated sensing, control, and actuation. At the core of this implementation lies the System Under Control (SUC), which encompasses key power-handling components: the battery switching mechanism, grid connection interface, and a DC/AC inverter responsible for converting stored energy into usable AC power. This configuration enables dynamic energy routing based on real-time demand and pricing conditions. The physical layout and interconnection of these elements are illustrated in Fig.~\ref{fig:prototype}. The subsequent subsections detail the mechanisms and control strategies used to manage these com

\subsection{Sensing Layer Description}

The sensing layer serves as the foundation of the energy management system by continuously monitoring electrical parameters that inform control decisions. At the core of its logic is the active power calculation:

\begin{equation}
    P = I \cdot V \cdot \text{PF}
\end{equation}

where $P$ denotes the active power, $V$ is the RMS voltage, $I$ is the RMS current, and $\text{PF}$ represents the power factor. Although this formula is elementary, it underpins the system's real-time monitoring capabilities. 
Voltage and current measurements are acquired using PZEM-016 modules. These modules compute not only active power but also reactive power, apparent power, and power factor. On the supply side, the battery subsystem is equipped with a dedicated BMS, which ensures both operational safety and optimal performance.
The energy storage component consists of a lithium iron phosphate (LiFePO$_4$) battery bank configured as 16 cells in series (16S), each with a nominal voltage of 3.2\,V, resulting in a total nominal pack voltage of 51.2\,V~\cite{battery}. The BMS monitors individual cell voltages, overall pack current, and temperature. It also estimates key metrics such as the SoC and State of Health (SoH), while enforcing critical protections against over-voltage, under-voltage, over-current, short circuits, and thermal events. When a fault is detected, the BMS immediately flags an alarm to the Raspberry Pi, enabling real-time response and user notification.

To maintain operational continuity during grid outages, the system incorporates an Uninterruptible Power Supply (UPS). This UPS supports the Raspberry Pi controller and critical relays, enabling an autonomous transition to battery power within approximately two seconds following a grid failure. With a total usable capacity of 4608\,Wh, the battery bank can sustain backup operation for durations dependent on load conditions and energy reserves.
The sensing layer collects data from both the supply side (grid and battery) and the demand side (load). These real-time measurements are transmitted to the Raspberry Pi 5 via RS-485 communication, ensuring robust and high-speed data transfer.
Using the collected data, the system enforces a predefined power threshold of 700\,W. When the measured load exceeds this limit, the control algorithm triggers a switch from grid power to battery storage, thereby implementing peak load mitigation strategies consistent with demand response and time-of-use (TOU) pricing models.

\subsection{Control Layer}

The control layer of the system is centered around a Raspberry Pi 5, a compact and versatile single-board computer widely adopted for embedded systems, automation, and IoT applications~\cite{pi}. In this project, the Raspberry Pi serves as the central controller, running custom-developed firmware written in Python~\cite{python}. Communication with the Battery Management System (BMS) is established via the Controller Area Network (CAN) protocol~\cite{can}, enabling real-time access to battery parameters such a SoC, voltage, current, and temperature.

Custom Python scripts process incoming CAN data, assess system status, and autonomously control relays to manage battery charging and discharging, based on predefined thresholds and electricity pricing schemes. This seamless software-hardware integration enables dynamic and adaptive operation under changing energy demands.
The control logic is implemented through six Python scripts and two MQTT brokers: one locally hosted and one cloud-based via the Blynk platform. The local broker handles real-time data acquisition from the BMS and two PZEM-016 sensors that monitor load and grid parameters. This configuration supports continuous operation even in the absence of internet connectivity. The Blynk broker facilitates bidirectional communication with the Blynk mobile application for real-time data visualization and remote control.

Two scripts are dedicated to acquiring and publishing BMS and PZEM data to the local MQTT broker. Two additional scripts forward these data streams to the Blynk cloud platform. The BMS script, based on Pylontech documentation~\cite{plewkaPylontech}, uses a \texttt{.dbc} file to decode CAN messages, while the PZEM script parses voltage, current, and power data. The Blynk interface is organized into distinct templates for battery monitoring, power sensing, and relay control.

Relay management is governed by two scripts: one interprets real-time system parameters to automate energy source switching, while the other processes user override commands from the Blynk app. The control logic prioritizes manual inputs when issued but defaults to autonomous operation under normal conditions. An emergency stop feature is also included, implemented via a physical pushbutton interfaced through the Raspberry Pi GPIO. The main control script (\texttt{main\_script\_spring25.py}) continuously monitors this input, immediately de-energizing all relays if triggered, and enabling safe recovery once cleared.

The core control logic is implemented in the \texttt{control\_relays()} function, which determines switching behavior based on SoC, load power (\texttt{pw\_load}), relay mode (\texttt{relay\_mode}), emergency status (\texttt{em\_mode}), and current time (\texttt{present\_minute}). Relay state transitions are handled by the functions \texttt{relays\_on()} and \texttt{relays\_off()}, with a 0.25-second delay preceding external DC relay activation to reduce arcing risks.

The decision-making process is formalized in Algorithm~\ref{alg:relay_control}, which prioritizes battery discharge during peak loads, charging during low SoC and off-peak periods, and enforces emergency shutdown when required.

\begin{algorithm}
\caption{Relay Control Logic}
\label{alg:relay_control}
\begin{algorithmic}[1]
\Require SoC, Present Minute (\texttt{PresMin}), Relay Mode (\texttt{RelMode}), Emergency Mode (\texttt{EmMode}), Load Power (\texttt{PwLoad}), Power Threshold (\texttt{PwThresh})
\Ensure Determine whether to enable battery mode (relays ON) or grid mode (relays OFF)

\If{\texttt{EmMode} == 0}
    \If{\texttt{RelMode} $\neq$ 1 and \texttt{RelMode} $\neq$ 2} \Comment{Auto Mode}
        \If{\texttt{PwLoad} < \texttt{PwThresh}}
            \If{\texttt{SoC} > 20\% and \texttt{PresMin} > 6}
                \State Set relays ON (battery mode)
            \ElsIf{\texttt{SoC} $\neq$ 100\% and \texttt{PresMin} < 6}
                \State Set relays OFF (grid + battery charging)
            \ElsIf{\texttt{SoC} == 100\% and \texttt{PresMin} < 6}
                \State Set relays ON (battery mode)
            \ElsIf{\texttt{SoC} < 20\%}
                \State Set relays OFF (grid charging)
            \EndIf
        \Else
            \State Set relays ON (battery mode)
        \EndIf
    \ElsIf{\texttt{RelMode} == 1}
        \State Set relays ON (force battery mode)
    \ElsIf{\texttt{RelMode} == 2}
        \State Set relays OFF (force grid mode)
    \EndIf
\Else
    \State Set relays OFF (emergency shutdown)
\EndIf

\end{algorithmic}
\end{algorithm}

In automatic mode, the controller uses a 700\,W threshold to determine switching. If the load exceeds this limit, the system shifts to battery power to minimize demand charges. Otherwise, it evaluates SoC and time-of-day metrics to determine the optimal mode, mimicking demand response strategies.

A graphical user interface (GUI), built using the \texttt{tkinter} library, provides real-time visualization of key parameters via the Raspberry Pi touchscreen. The GUI is updated by three primary functions: \texttt{update\_readings()}, \texttt{update\_font\_size()}, and \texttt{exit\_fullscreen()}. These functions ensure continuous display of SoC, load power, grid status, and relay states, offering users an intuitive, real-time interface for system monitoring.

\subsection{Actuator Layer}

The actuator layer is responsible for executing control commands by physically regulating the flow of electrical power between the grid, battery, inverter, and system loads. It ensures seamless switching between power sources, supports safe operation, and enforces charging or discharging behavior based on system logic.

The power path begins with AC input from the utility grid. The neutral line is routed to the AC common block, the ground line to the AC ground block, and the hot line passes through the normally closed terminal of an emergency stop button. Under normal conditions, this terminal remains closed, allowing uninterrupted power delivery. In the event of an emergency, pressing the button immediately disconnects the AC input, ensuring user and system safety.

Beyond the emergency stop, the hot conductor feeds into the main AC voltage block, which distributes power throughout the system. Source selection between the grid and the inverter is managed by Relay 4 and Relay 5, both situated on the PCB relay board. Specifically, Relay 4 controls power flow to the system loads. Its common terminal connects to the loads, the normally closed terminal connects to the grid input, and the normally open terminal connects to the inverter output. Relay 5 manages the UPS input using a similar configuration.

Relay 6 governs the battery's operational mode by toggling between charging and discharging states. Due to the higher current requirements in this path, Relay 6 actuates a large external relay. This external relay, in turn, routes battery output either to the inverter for discharging or to the battery charger for recharging, depending on current operating conditions. The control voltage for both Relay 6 and the external relay is supplied by a 12-volt voltage block, which shares a common ground reference with the DC ground block to ensure electrical stability.
The 12-volt block also supplies power to the PCB relay board and the Raspberry Pi, maintaining low-voltage control system operation. During power outages, the UPS provides backup 12-volt power, preserving relay actuation and control logic functionality.

System grounding is implemented through the AC ground block, while a stable neutral reference across the UPS, inverter, battery charger, and system loads is maintained via the AC common block. Battery wiring is designed for safety, with the positive terminal routed through a fuse before reaching the external relay's common terminal. The negative terminal is securely bonded to the DC ground block. Proper coordination between the inverter and battery charger ensures consistent power flow during both grid-connected and backup operating modes. This actuator architecture enables safe and reliable switching, maintains stable voltage references, and adheres to standard electrical safety practices across all operating scenarios.

\section{Experimental Results}

The proposed system underwent comprehensive experimental validation to assess the accuracy of real-time measurements, evaluate relay switching behavior, and verify overall system performance under realistic operating conditions.

\subsubsection{Measurement Accuracy Validation}

To evaluate the accuracy of the PZEM-016 energy monitoring modules, known electrical loads were connected to the system, and their measured values were compared against expected outputs. For instance, a 7\,W light bulb produced a corresponding reading of approximately 7\,W at the Raspberry Pi interface, confirming the modules' ability to accurately capture voltage and current data. 

According to the PZEM-016 specifications~\cite{pzem_016_amazon}, measurement deviations were within the expected resolution range of approximately $\pm 0.1$\,W. Additionally, polling interval tests revealed that the minimum reliable sampling period for the PZEM-016 is one second. Attempts to reduce the interval to 0.5 seconds resulted in repeated measurements, indicating a limitation in the sensor’s internal refresh rate.

\subsubsection{Relay Switching and System Behavior}

\begin{figure}[htbp]
    \centering
    \includegraphics[width=0.8\linewidth]{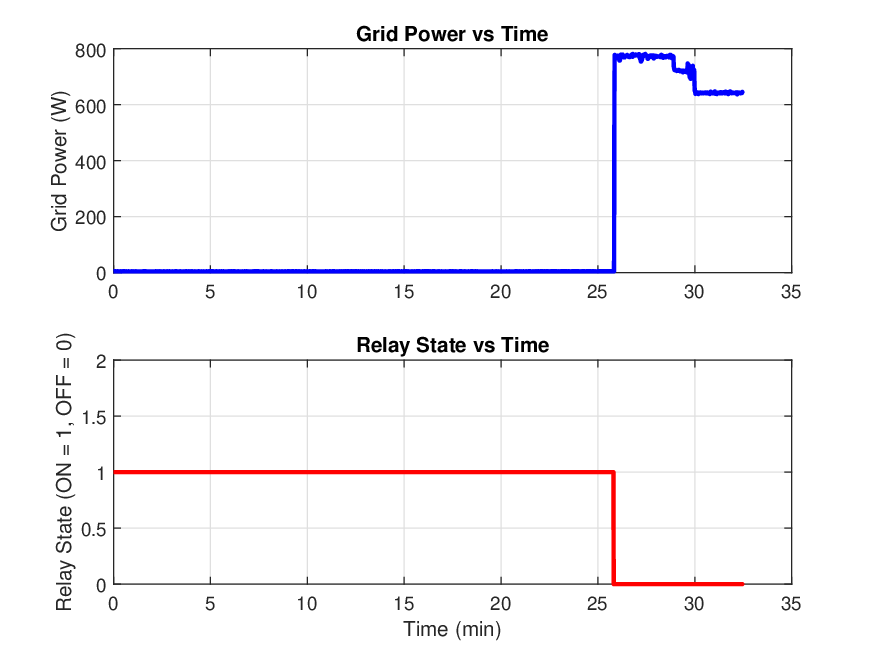}
    \caption{System switching from battery power to grid power}
    \label{fig:Grid_to_Batt_Mode}
\end{figure}

Relay switching behavior was evaluated by monitoring system response during mode transitions. As illustrated in Figure~\ref{fig:Grid_to_Batt_Mode}, the system initially operated in battery mode and then transitioned to grid supply upon receiving a remote command via the Blynk application. This transition resulted in a distinct increase in grid power consumption, attributed to both active load support and concurrent battery charging operations.

During battery mode, a small residual grid power draw of approximately 4.7\,W was observed. This behavior is consistent with standby current requirements from system components such as the inverter, battery charger, and PZEM modules, aligning with previously reported standby characteristics~\cite{pzem_016_amazon}.

\subsubsection{State of Charge Monitoring}

\begin{figure}[htbp]
    \centering
    \includegraphics[width=1\linewidth]{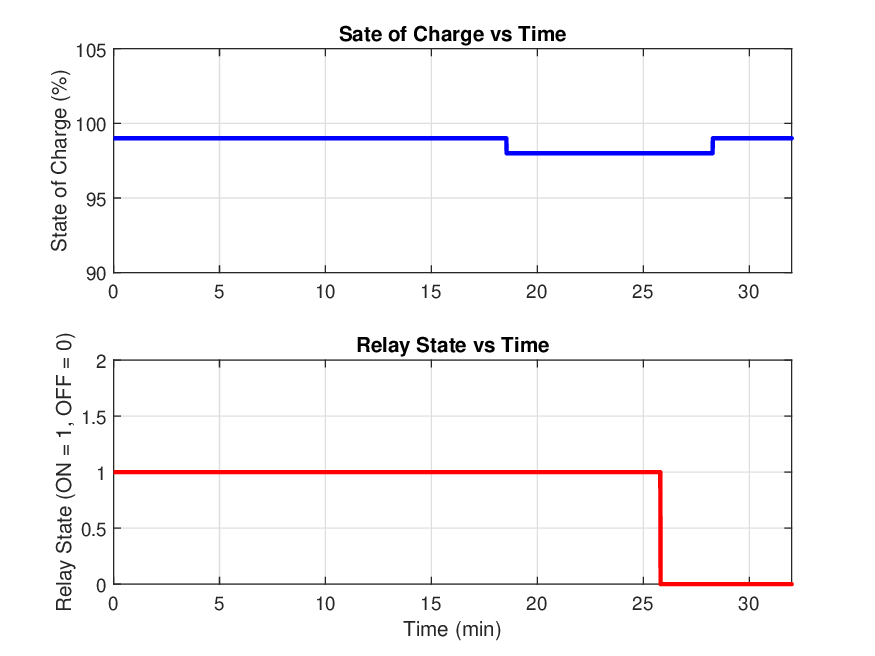}
    \caption{Battery mode to grid mode system with battery discharging and charging}
    \label{fig:SOC_Relay_Status}
\end{figure}

State of Charge monitoring functionality was validated by comparing the SoC values reported on the battery's built-in display with those displayed on the Raspberry Pi interface. As illustrated in Fig.~\ref{fig:SOC_Relay_Status}, the SoC decreased during battery discharge operation and increased once the system transitioned back to grid power and initiated charging. This behavior confirms accurate real-time integration between the BMS and the control interface, ensuring reliable monitoring and control of battery energy levels.

\subsubsection{Power Factor Performance Analysis}

\begin{figure}[htbp]
    \centering
    \includegraphics[width=0.9\linewidth]{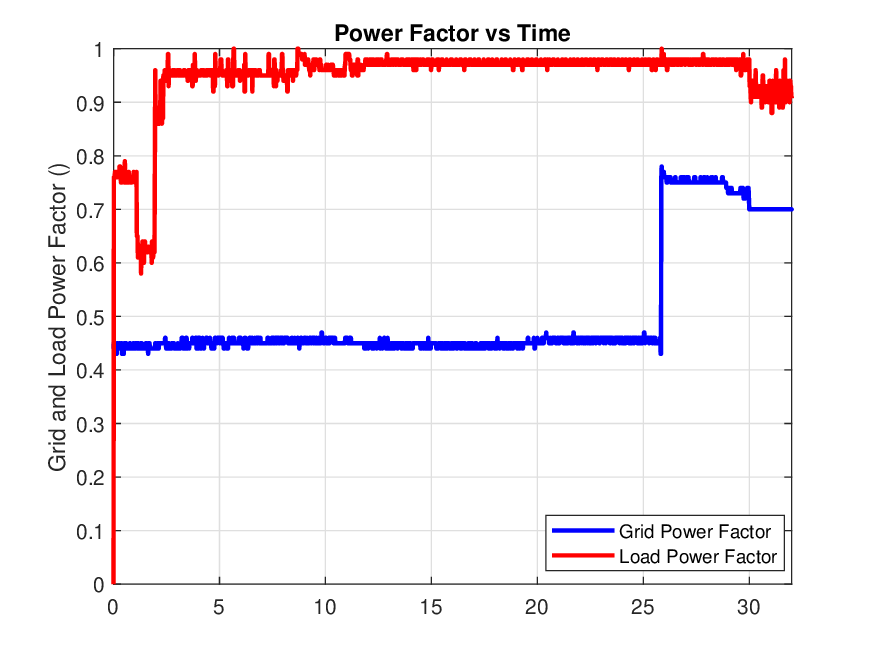}
    \caption{Grid and load power factor under battery and grid operation}
    \label{fig:Load_Grid_PF}
\end{figure}

System power factor behavior under different operating modes is illustrated in Figure~\ref{fig:Load_Grid_PF}. The load-side power factor remained consistently high, ranging between 0.95 and 0.98, during both battery and grid operations, indicating efficient utilization of apparent power at the load level.
In contrast, the grid-side power factor was initially low, approximately 0.45, during battery mode. This reduction is attributed to minimal active power draw combined with the presence of inductive standby components such as the inverter and battery charger. When the system transitioned to grid mode, the grid-side power factor increased to around 0.75, reflecting higher real power demand and improved reactive power balance.

These observations suggest opportunities for future optimization, particularly through the electrical isolation of standby devices like the inverter and charger during pure battery operation. Such modifications could help minimize reactive power contribution and enhance overall system efficiency.

\subsection{Performance}

The proposed system employs a stacked logic flow architecture for decision-making, leveraging real-time data streams to enable reliable, adaptive power distribution. Unlike conventional systems that rely on static timers or manual switching, the implemented eco-mode autonomously evaluates live power measurements from PZEM sensors in conjunction with battery State of Charge (SoC) and power factor data, all communicated via MQTT protocol.

\begin{figure}[htbp]
    \centering
    \includegraphics[width=0.95\linewidth]{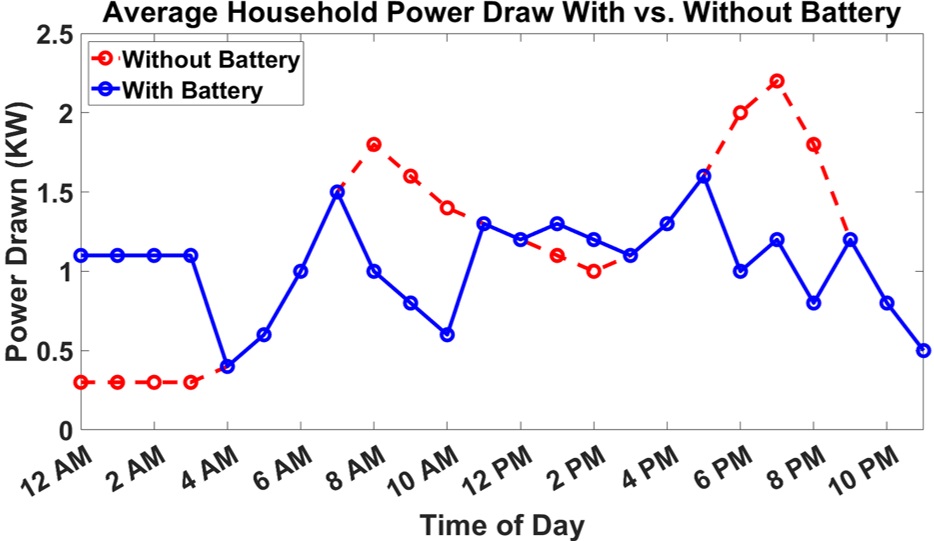}
    \caption{Average household power draw with and without battery support.}
    \label{fig:Vs}
\end{figure}

Experimental validation was performed using a 5120\,Wh lithium iron phosphate (LiFePO$_4$) battery, which demonstrated substantial backup capacity for residential energy applications. For an average household load of approximately 1300\,W, the battery was capable of sustaining critical operations for roughly 3.9 hours during an outage scenario.
Further insights were obtained through a MATLAB-based simulation that modeled battery discharge and recharge cycles aligned with utility pricing patterns. The system discharged 5.4\,kWh during high-demand (peak tariff) windows and absorbed 3.6\,kWh during low-demand (off-peak) periods. As shown in Figure~\ref{fig:Vs}, coordinated battery scheduling effectively flattened the household load curve, thereby mitigating grid stress and enhancing energy cost efficiency. In addition, a 37\,Wh UPS was integrated to support the Raspberry Pi 5 controller and essential relay control circuitry during full system outages. With an average power draw of approximately 8.5\,W, the UPS maintains system operability for approximately 4 hours and 20 minutes. This enhances system fault tolerance and ensures safe, autonomous shutdown when extended outages occur.

\section{Conclusion}
\label{sec:conslusion}

This paper has presented the development, implementation, and experimental validation of a smart, low-cost Energy Management System (EMS) and Demand Charge Management (DCM) prototype. Designed as a practical and scalable alternative to expensive commercial systems, the proposed solution integrates real-time sensing, Python-based control logic, MQTT communication, and IoT visualization through the Blynk platform. It dynamically manages energy usage by switching between grid and battery power based on real-time measurements of load demand, electricity pricing, and battery State of Charge (SoC), thereby achieving both economic optimization and operational resilience.
Experimental results confirmed the system's ability to reduce electricity costs through peak demand mitigation and energy arbitrage, while also providing reliable backup power during grid outages. Performance testing validated the accuracy of power measurements, responsiveness of relay switching, and the robustness of the control logic. A MATLAB-based simulation further demonstrated the system's effectiveness in flattening residential load profiles, reducing grid stress, and enhancing energy efficiency.
In addition to its technical capabilities, the system is designed to serve as an educational platform for undergraduate instruction in power systems, embedded control, and energy storage. Developed through a senior design project, the system exemplifies the application of interdisciplinary engineering knowledge in a real-world context. Its affordability, modular design, and use of open-source tools make it especially well-suited for adoption in teaching laboratories and student-driven research.
The key contributions of this work include the design of a complete, three-layer EMS/DCM hardware architecture; the integration of low-cost, open-source technologies for real-time control and monitoring; and the successful validation of system functionality under realistic operating conditions. Moreover, the project illustrates how hands-on prototyping can bridge theoretical knowledge with practical implementation, offering value not only for distributed energy management but also for engineering education.
Future enhancements may include scaling the system to support multi-node microgrids, integrating renewable energy sources such as solar photovoltaics, and incorporating machine learning algorithms for predictive energy scheduling and adaptive control.

\bibliographystyle{IEEEtran}
\IEEEtriggeratref{50}
\bibliography{ref}
\end{document}